\documentclass[12pt]{article}
\usepackage{amsmath,amssymb}

\newcommand{\un}{\mathbf{1}}
\newcommand{\R}{\mathbb{R}}

\newcommand{\eps}{\varepsilon}

\newcommand{\dt}{\partial_t}

\newtheorem{theo}{Theorem}
\newtheorem{prop}[theo]{Proposition}
\newtheorem{lemm}[theo]{Lemma}
\newtheorem{corr}[theo]{Corollary}

\newtheorem{defi}[theo]{Definition}
\def\qed{\hbox{${\vcenter{\vbox{
  \hrule height 0.4pt\hbox{\vrule width 0.4pt height 6pt
  \kern5pt\vrule width 0.4pt}\hrule height 0.4pt}}}$}}
\def\ds{\displaystyle}

\def\begproof{\noindent {\bf Proof. }}
\def\blackbox{\leavevmode\vrule height 5pt width 4pt depth 0pt\relax}
\def\endproof{\null\hfill {$\blackbox$}\bigskip}

\long\def\symbolfootnote[#1]#2{\begingroup%
\def\thefootnote{\fnsymbol{footnote}}\footnote[#1]{#2}\endgroup}
\begin{document}


\begin{center}
{\Large\sc Global regularity of solutions  to systems of reaction-diffusion  with Sub-Quadratic Growth
in any dimension}
\\
\vspace*{20pt} M. Cristina {\sc Caputo}$^1 ${}, Alexis {\sc
Vasseur}$^{2}$\symbolfootnote[2]{A. Vasseur was partially supported
by the NSF, DMS- 0607053.}

\vspace*{15pt} $^{1}$
Department of Mathematics\\
University of Texas at Austin\\
Austin, Texas 78712,
USA\\
E-mail: {\tt caputo@math.utexas.edu}
\\
\vspace*{15pt} $^{2}$
Department of Mathematics\\
University of Texas at Austin\\
Austin, Texas 78712,
USA\\
E-mail: {\tt vasseur@math.utexas.edu}
\end{center}

\vspace{0.5 cm}

\begin{abstract}
This paper is devoted to the study of the regularity of solutions to
some systems of reaction--diffusion equations, with reaction terms having a subquadratic growth. We
show the global  boundedness and regularity of  solutions, without smallness assumptions, in any dimension $N$.
 The proof is based on blow-up techniques.
The
natural entropy of the system plays a crucial role in the analysis.
It allows us to use of De Giorgi type methods introduced for elliptic regularity with rough coefficients.
In spite these systems are entropy supercritical, it is possible to control the hypothetical blow-ups, in the critical scaling, via
a very weak norm. 
Analogies with the Navier-Stokes equation are briefly discussed in the introduction.
\end{abstract}

\medskip
\noindent {\bf Key words.}
Reaction-diffusion systems. Global regularity. Blow-up methods.
 \vskip 0.4cm

\medskip
\noindent {\bf AMS Subject classification.}
 35Q99, 
 35B25, 
 82C70. 
\vskip 0.4cm


\section{Introduction}

This paper is dedicated to the study of the global regularity in time of the solutions to a class of reaction-diffusion systems. Reaction-diffusion systems are used as models for a variety of problems, especially  in chemistry and  biology \cite{Jef, ErTo, Feng, Fife, Mo2, Mu}. The question of the existence of global solutions is particularly important and it has been widely studied \cite{Mo1, Mo4, Mo5, Rot, SoPe}.
In full generality, such system can have solutions which blow up in finite time \cite{PiSc} as it is well known when considering non linear heat equations  \cite{GiKo, We}. In this article we focus on systems which are entropy supercritical, that is, such that the preserved physical quantities, as the mass and the entropy, are not shrinking (not even preserved) via the universal scaling of the system. The global regularity for entropy (or energy) supercritical problems is an important question in several areas like nonlinear waves, Schroedinger, and the Navier-Stokes equations in dimension $N\geq 3$. 
\vskip 0.2cm
We present, here, a class of systems for which we are able to show global regularity despite their supercritical nature.

We consider systems of the
form:
\vskip 0.2cm
\begin{equation}\label{eq_system}
\left\{\begin{array}{l}
\displaystyle{\dt a-D\Delta a=Q(a), \qquad t>0, \ x\in
\R^N,}\\[0.3cm]
\displaystyle{a(0,x)=(a_1^0(x), a_2^0(x),\cdot\cdot\cdot, a_P^0(x)),
\qquad x\in \R^N,\qquad P\geq 1.}
\end{array}
\right.
\end{equation}
\vskip 0.4cm
\noindent We assume that the matrix of diffusion is diagonal:
 $$
D=\mathrm{diag}(D_1,\cdot\cdot\cdot,D_P),\qquad D_i\in\R, \hbox{ for $\,1\leq i\leq P$},
 $$
 with
 $$
\underbar d:=\inf_{1\leq i\leq P}\{D_i\}>0\qquad \hbox{and}\qquad  \bar d:=\sup_{1\leq i\leq P}\{D_i\}.
 $$
 In addition, we assume that the reaction diffusion term $Q$ is regular and  it satisfies the following four conditions:
there exists $0<\nu<2$ and $\Lambda>0$ such that

\begin{align}
 & Q_i(a)\geq 0 ,\,\,\,\,\hbox{if}\,\,\,a_i\leq0,\,\,\,\,  \,1\leq i\leq P,&\,\,\,\, \qquad \qquad \qquad\label{Hypothesis_positivity}\\
& |\nabla Q_i(a)| \leq \Lambda |a|^{\nu-1}, \qquad a\in\R^P_+,\,\,\,\,  \,1\leq i\leq P, &\,\,\,\,\qquad \qquad \label{Hypothesis_subquadratic}\\
& \sum_{i=1}^{P}Q_i(a) =0,  \qquad     a\in
\R^P_+, &\,\,\,\,\qquad\qquad \qquad  \label{Hypothesis_mass} \\
& \sum_{i=1}^{P}\ln a_iQ_i(a)\leq0, \qquad   a\in \R^P_+.&\,\,\,\,\qquad \qquad \qquad \label{Hypothesis_entropy}
\end{align}

Hypothesis (\ref{Hypothesis_positivity}) ensures the nonnegativity
of the $a_i$,  hypothesis
(\ref{Hypothesis_subquadratic}) is a restriction on the growth of
$Q$,
hypothesis (\ref{Hypothesis_mass}) ensures the
conservation of the total mass $\sum_{i=1}^{P}\int a_i\,dx$, and
hypothesis (\ref{Hypothesis_entropy}) ensures the non increase of
the entropy $\sum_i\int  a_i\ln a_i\,dx$. Note that (\ref{Hypothesis_subquadratic})
is the only restrictive hypothesis. It requires that the system is at most subquadratic. The other hypothesis are physical and pretty standard in systems coming from chemistry.\\

Examples of systems verifying such hypothesis for $P=4$ are given by:

$$Q_i(a)=(-1)^i(\phi(a_1a_3)-\phi(a_2a_4)),\qquad a=(a_1,a_2,a_3,a_4)\in \R^4_+, \,\,\,1\leq i\leq 4,$$

where $\phi\in C^\infty(\R^+)$ is any nondecreasing function verifying
\begin{eqnarray*}
&& \phi(z)=0,\qquad \mathrm{for}\  \ z<0,\\
&&\phi(z)=z^{\nu/2},\qquad \mathrm{for}\  \ z>1.
\end{eqnarray*}

\vskip1cm
We consider nonnegative initial  values   $a_i^0\in C^{\infty}(\R^N)\cap L^{\infty}(\R^N)$, $1\leq i\leq P$,  verifying
\begin{equation}\label{M_00}
\sup_{1\leq i\leq P} \left\{\int_{\R^N}\,a^0_i(1+|x|+|\ln(a^0_i)|)\,dx\right\}<\infty,
\end{equation}
but without smallness condition. We note that these additional constraints correspond to the finite initial mass and total entropy, together with a condition of confinement of the mass near the origin. \\

The main result of this paper is the following theorem.

\begin{theo}\label{Theo_principal}
Let $\nu<2$. Consider a system~(\ref{eq_system}) verifying
hypothesis~
(\ref{Hypothesis_positivity}), (\ref{Hypothesis_subquadratic}),
 (\ref{Hypothesis_mass}), and
 (\ref{Hypothesis_entropy}).
Then, for any  $a^0=(a_1^0,\cdot\cdot\cdot, a_P^0)$, with $a_i^0\in C^{\infty}(\R^N)\cap L^{\infty}(\R^N)$, $a_i^0\geq 0$ for $\,1\leq i\leq P$, verifying the condition~(\ref{M_00}),
there
exists a unique smooth global solution $a=(a_1,\cdot\cdot\cdot, a_P)$ defined on $[0,\infty)\times\R^N$.
Moreover, $a_i\geq 0$  for all  $\,1\leq i\leq P$.
\end{theo}
A system of particular interest, which is a limit case (not included) of this study, is the quadratic  system of reaction with four species corresponding to the chemical reaction
$$
A_1+A_3\rightleftharpoons A_2+A_4.
$$
 In this case the system~(\ref{eq_system}) has
$P=4$, and
$$
Q_i(a)=(-1)^i(a_1a_3-a_2a_4),\qquad a=(a_1,a_2,a_3,a_4)\in \R^4_+, \,\,\,1\leq i\leq 4,
$$
where $a_i$ is the concentration of the species $A_i$.
Let us denote
$$
U_i=\sqrt{a_i},\qquad 1\leq i\leq 4.
$$
The mass conservation and entropy dissipation give (see section \ref{sec:entropy}) the following bounds for any $T>0$:
$$
U\in L^\infty(0,T;L^2(\R^N)),\qquad \nabla_x U\in L^2((0,T)\times\R^N).
$$
Also, the universal scaling is given by
$$
U_\eps(t,x)=\eps U(\eps^2t,\eps x), \qquad \eps>0.
$$
This means that $U$ is a solution to the system~(\ref{eq_system}) if and only if $U_\eps$ is a solution to the same problem.
Note that $U$ has exactly the same conserved quantities  and the exact same universal scaling as the Navier-Stokes equation, which is supercritical for $N>2$.
We interpret this by saying that our result gives the full regularity for a family of parabolic systems, in any dimension, which has almost the same supercriticality as the Navier-Stokes equation.

 However, we observe that these systems of reaction-diffusion have some features quite different from fluid mechanics. In the inviscid case where $D_i=0$ for all $i$,  $\,1\leq i\leq P$,  a maximum principle  ensures that the solutions remain globally bounded. But this structure breaks down in the diffusive case provided that the diffusion coefficients $D_i$'s are not all the same. Surprisingly, such systems are, actually, destabilized by the diffusion.
\vskip0.1cm

 One of the key arguments in our proof is the use of a local parabolic regularization effect. This argument basically claims that if the nonlinear terms are somehow  controlled and small in a cylinder $(t_0-1,t_0)\times B(x_0,1)$, then the value of the solution at $(t_0,x_0)$ is bounded by 1. This argument is at the heart of the partial regularity results for the Navier-Stokes equations \cite{CKN, Lin, Sch, Sch1}.
In \cite{Va}, a new proof is proposed based on the De Giorgi's parabolic regularity method \cite{DeG}. This technique is particularly powerful and gives important results in different physical areas, such as the quasigeostrophic equation \cite{CaVa}, and the compressible Fourier-Navier-Stokes system \cite{MeVa1, MeVa2}. This method has been used for the first time in the context of reaction-diffusion equation in \cite{GoVa}. This paper contains the cases $N\leq 2$, even with quadratic growth. It corresponds to the entropy critical case.

The local parabolic regularization results are based on recursive controls, along a family of shrinking cylinders, of the nonlinear terms in a cylinder by the entropy (or energy) on a bigger one. We note that, in the context of reaction-diffusion systems, the classical methods based on the Green function of the heat equation do not work in the supercritical cases. The main problem is that we have different diffusion coefficients. In this context, the De Giorgi method proves to be particularly powerful as it exploits the physical entropy quantity~(\ref{Hypothesis_entropy}). Using this property, this method depletes the nonlinearity of one exponent. This is one of the key facts which allows us to work with subquadratic nonlinearity for any dimension $N$. It follows that, locally, it is enough to control recursively slightly bigger (to the $\mathrm{log}$) norms that the $L^1$ norm to get the regularization since such norms are controlled by the entropy.

\vskip0.1cm

To get the global regularity, we shrink these local norms, about any point $(t_0,x_0)$ through the universal scaling
\[
a^{\eps}(s,y)=\eps^{\frac{2}{\nu-1}}a(\eps^2 s+t_0,\eps y+x_0), \,\,\,\,\hbox{for}\,\,\, \eps>0.
\]
If $a=(a_1,\cdot\cdot\cdot, a_P)$ is a solution to the system, then $a^\eps=(a_1^\eps,\cdot\cdot\cdot, a_P^\eps)$ is a solution to a system with a reaction term $Q^\eps$ which shares the same properties as the reaction term $Q$ [Hypothesis~(\ref{Hypothesis_positivity}), (\ref{Hypothesis_subquadratic}),  (\ref{Hypothesis_mass}), and (\ref{Hypothesis_entropy})] with the same $\nu$ and the same constant $\Lambda$.
If there is an $\eps>0$, small enough,  for which the norm is very small, then we can use the local parabolic regularization result on $a^\eps=(a_1^\eps,\cdot\cdot\cdot, a_P^\eps)$ to ensure the regularity of $a=(a_1,\cdot\cdot\cdot, a_P)$ at the point $(t_0,x_0)$.
In the entropy supercritical cases, all the quantities obtained through the entropy blow up along the rescaling. But, on the other hand, we can show than the solution is globally bounded in a very weak space which has the same homogeneity as $L^\infty(W^{-2,\infty})$. The norm in this weak space shrinks in the subquadratic growth $\nu<2$ as
$$
\|a^\eps\|_{L^\infty(W^{-2,\infty})}=\eps^{\frac{4-2\nu}{\nu-1}}\|a\|_{L^\infty(W^{-2,\infty})}.
$$
It is the nonnegativity of the solutions that allows us to control the depleted nonlinear term by this weak norm. We can use, then,  the local parabolic regularization property.
\vskip0.1cm

The control of the solution in this weak space is based on a nice duality argument introduced by Pierre and Schmitt~\cite{PiSc}. In~\cite{PiSc}, they obtain a global bound in $L^2(L^2)$ of the solution. A similar argument has been used in~\cite{DFPV} to get global weak solutions. The $L^2(L^2)$ norm does not shrink through the universal scaling for the  entropy supercritical cases. However, we 
will show that the same method can be applied for the weak norm that we present. 

\vskip0.1cm

The study of this weak norm was suggested by the paper \cite{Va2}, where the $L^\infty(W^{-1,\infty})$ norm is shown to play an important role in the Navier-Stokes equation. We observe that  the $L^\infty(W^{-1,\infty})$ norm has the same scaling as the $L^\infty(BMO^{-1})$ norm studied by Koch and Tataru for the Navier-Stokes equation. The weak norm that we introduce corresponds, formally, for the quadratic case, to the $L^\infty(W^{-2,\infty})$ norm.

The additional control of the weak norm comes from a linear equation verified by the total mass $\rho=\sum_i^P a_i$. Pierre and Schmitt~\cite{PiSc} use this extra equation to get the $L^2(L^2)$ estimates. They also show that the full regularity cannot come solely from this argument as they provide some explicit examples of solutions to this type of linear equation which blow up in finite time.

\vskip0.1cm

We note that, in the contest of the Navier-Stokes equation, 
an equation on the vorticity  is also provided as an extra. For the Euler equation (the inviscid case), the vorticity $\omega$ is a solution to the equation:
$$
\dt \omega+(u\cdot\nabla) \omega-(\omega\cdot\nabla)u=0.
$$
If we drop the dependence of vorticity on the velocity $u$, the above equation can be seen as a linear equation on $\omega$, which can be solved in the Lagrangian coordinates as
$$
\omega(t,a)=(\omega(0,a)\cdot\nabla_a)X(t,a),
$$
where $X$ is the flow given by
$$
\dt X(t,a)=u(t,X(t,a)).
$$
Note that, for any smooth initial values, this gives a uniform control of the vorticity in $L^\infty(W^{-1,\infty}_{\mathrm{loc}})$ in the Lagrangian coordinates. This norm is, actually, shrinking through the universal scaling of the Navier-Stokes equation. Unfortunately, this structure seems to be destroyed by the viscosity term. Such proprerty provides  another analogy with the reaction-diffusion systems: as for the maximum principle for the reaction-diffusion systems, a crucial supercritical structure known on the Euler system is destabilized by the diffusion term.
\vskip0.1cm

Global existence of weak solutions of \eqref{eq_system}, was established in~\cite{DFPV}.
The dissipation property \eqref{Hypothesis_entropy}
is also the basic tool for studying the asymptotic trend to equilibrium~\cite{DeFe2, DeFe1} in the spirit
of the entropy/entropy dissipation techniques which are presented e.g. in \cite{Vil}
(we refer also to \cite{Ball} for further investigation of the large time behavior of nonlinear evolution  systems using the entropy dissipation).
Let us also mention
that \eqref{eq_system} can be derived through hydrodynamic scaling from kinetic models, see \cite{BiDe}.


%

%

\vskip1cm
As remarked above, the results of Theorem~\ref{Theo_principal} are trivial if the diffusion coefficients $D_i$'s are all equal. Another trivial case corresponds to $P=2$ (two species) where a maximum principle holds (even without the subquadratic property (\ref{Hypothesis_subquadratic})).  For the sake of completeness we will give a proof of this result in the appendix as we did not find this case in the literature.
\vskip0.4cm
\begin{theo}\label{theo_max_principle}
Consider a system~(\ref{eq_system}), with $P=2$, and  verifying
hypothesis~
(\ref{Hypothesis_positivity})
 (\ref{Hypothesis_mass}), and
 (\ref{Hypothesis_entropy}).
Then, for any  $a^0=(a_1^0,a^0_2)$, with $a_1^0,a_2^0\in C^{\infty}(\R^N)\cap L^{\infty}(\R^N)$, $a_1^0,a_2^0\geq 0$,  verifying the condition~(\ref{M_00}),
there
exists a unique smooth global solution $a=(a_1,a_2)$ defined on $[0,\infty)\times\R^N$.
Moreover, $a_1,a_2\geq 0$, and
$$
\sup_{(t,x)\in(0,\infty)\times\R^N}\{a_1(t,x),a_2(t,x)\}\leq \sup_{x\in\R^N}\{a^0_1(x),a^0_2(x)\}.
$$
\end{theo}
\vskip0.2cm
The proof of the maximum principle collapses for $P\geq 3$.

\vskip1cm

Some of the estimates in this paper will be based on results obtained in~\cite{GoVa}, but,  for the sake of completeness this paper is self contained.

For any smooth initial data,
 standard theory ensures the existence of a smooth solution on (at least) a short time $[0,T)$, $T>0$. We set $T_0$ to be the biggest of such lapse of time.
 Our aim is  to show that  $T_0=\infty$. Standard bootstrapping arguments give  that if $T_0<\infty$ then
 $$
 \lim_{t\to T_0}\sup_{1\leq i\leq P}\|a_i(t,\cdot)\|_{L^\infty(\R^N)}=+\infty.
 $$
 We will obtain a uniform bound on $[0,T_0)$ contradicting the blow-up of the solution in finite time. Indeed,
 we will show that this bound depends  only on $T_0$ and the quantity $M_0$ defined by
\begin{equation}\label{M_0}
M_0=\sup_{1\leq i\leq P} \left\{\int_{\R^N}\,a^0_i(1+|x|+|\ln(a^0_i)|)\,dx+\|a^0_i\|_{L^\infty(\R^N)}\right\}.
\end{equation}

In the second section, we provide standard a priori estimates based on the mass conservation and the entropy dissipation. The third section is dedicated to the local parabolic regularization principle. The duality arguments are given in the fourth section. We introduce and then apply the rescaling arguments in the last section.

\section{Entropy Dissipation}\label{sec:entropy}
\vskip1cm
In this section, we derive a priori estimates on $[0,T_0)$ where the solution is smooth. The dimensional cases of $N\leq 2$ have been already studied in~\cite{GoVa}, but they can also be deduced from the techniques that we present in this paper. (The reader can easily replace the classical inequalities in the sequel of the proofs for dimensions $N\leq 2$). Henceforth, we will analyze here only the dimensions $N\geq3$. 
\\

We discuss the a priori estimates that can be naturally deduced from \eqref{Hypothesis_mass} and \eqref{Hypothesis_entropy}.
\vskip0.1cm
\begin{prop}\label{P1}
There exist two constants $C_0$ and $C_1$, such that the following is true. Let
~(\ref{eq_system}) be any system verifying (\ref{Hypothesis_positivity}), (\ref{Hypothesis_mass}), and (\ref{Hypothesis_entropy}), and
any initial values 
 $a_i^0\geq 0$, smooth, verifying
\begin{equation}\label{hypdatap}
M_0=\ds\sup_{1\leq i\leq P}\left(\|a_i^0\|_{L^\infty(\R^N)}+\ds\int_{\R^N} a_i^0\big(1+|x|+|\ln(a_i^0)|\big)\, dx\right)<\infty.
\end{equation}
Let $a=(a_1,\cdot\cdot\cdot, a_P)$ be the associated solution. Assume that  $a$ is regular on its maximal lapse of time $[0,T_0)$.
Let 
\[
\mathfrak D(a)=-\sum_{i=1}^{P}Q_i(a)\ln(a_i)
\geq 0.
\]
Then,  for any $0<T<T_0$,
\[
\begin{array}{l}
\ds\sup_{0\leq t\leq  T}\Big\{\ds\sum_{i=1}^{P}\ds\int_{\R^N}
a_i\big(1+|x|+|\ln(a_i)|\big)(t,x) \, dx\Big\} +2\underbar{d}
\ds\sum_{i=1}^{P}\ds\int_0^T\int_{\R^N} \big|\nabla_x
\sqrt{a_i}\big|^2dx\,ds
\\
\qquad\qquad\qquad\qquad\qquad\qquad\qquad\qquad
+
\ds\sum_{i=1}^{P}\ds\int_0^T\ds\int_{\R^N}
\mathfrak D(a) dx\,ds

\leq (C_0+C_1T)(M_0+1).\end{array}
\]


\end{prop}
\begproof
\noindent First, thanks to (\ref{Hypothesis_positivity}), we have $a_i(t,x)\geq0$ for any $t\in [0,T_0)$, $x\in\R^N$, and $1\leq i\leq P$. As a consequence of
\eqref{Hypothesis_mass} and \eqref{Hypothesis_entropy}, we get
\[
\ds\frac{d}{dt}\ds\sum_{i=1}^{P}\ds\int_{{\R^N}} a_i\,dx=0
\]
and
\[
\ds\frac{d}{dt}\left(
\ds\sum_{i=1}^{P}\ds\int_{\R^N} a_i\big(1+\ln(a_i)\big) \,dx\right)
+
\ds\sum_{i=1}^{P}\ds\int_{\R^N} D_i\nabla_x a_i\cdot\ds\frac{\nabla_x a_i}{a_i} \,dx
+\ds\int_{\R^N} 
 \mathfrak D
 \,dx
=0.
\]
Also,
\[\ds\sum_{i=1}^{P}\ds\int_{\R^N} D_i\nabla_x a_i\cdot\ds\frac{\nabla_x a_i}{a_i} \,dx
\geq \underbar{d} \ds\sum_{i=1}^{P}\ds\int_{\R^N} \ds\frac{|\nabla_x a_i|^2}{a_i} \,dx
=4\underbar{d} \ds\sum_{i=1}^{P}\ds\int_{\R^N} \big|\nabla_x \sqrt{a_i}\big|^2 \,dx.
\]
Moreover
\[\begin{array}{lll}
\ds\frac{d}{dt}
\ds\sum_{i=1}^{P}\ds\int_{\R^N} a_i|x|\,dx
&=&-
\ds\sum_{i=1}^{P}\ds\int_{\R^N} D_i\nabla_x a_i\cdot\ds\frac{x}{|x|}\,dx\\
&\leq&
\bar{d}\ds\sum_{i=1}^{P}\ds\int_{\R^N} |\nabla_x a_i|\,dx
=
\bar{d}\ds\sum_{i=1}^{P}\ds\int_{\R^N} \ds\frac{|\nabla_x a_i|}{\sqrt{a_i}}\sqrt{a_i}\,d x\\
&\leq&
\ds\frac{\underbar{d}}{2}\ds\sum_{i=1}^{P}\ds\int_{\R^N} \ds\frac{|\nabla_x a_i|^2}{a_i} \,dx
+
\ds\frac{\bar{d}^2}{2\underbar{d}}\ds\sum_{i=1}^{P}\ds\int_{\R^N} a_i\,dx.
\end{array}\]
Let $0\leq t\leq T$. We integrate on the interval $[0,t]$,
\[\begin{array}{l}
\ds\sum_{i=1}^{P}\ds\int_{\R^N} a_i\big(1+|x|+\ln(a_i)\big) \,dx
+\ds\frac{\underbar{d}}{2}\ds\sum_{i=1}^{P}\ds\int_0^t\ds\int_{\R^N} \ds\frac{|\nabla_x a_i|^2}{a_i} \,dx\,ds
+
\ds\int_0^t\ds\int_{\R^N}
\mathfrak D
\,dx\,ds
\\
\qquad\qquad\leq
M_0 + \ds\frac{\bar{d}^2}{2\underbar{d}}\ds\sum_{i=1}^{P}\ds\int_0^t\ds\int_{\R^N} a_i\,dx\,ds
\\
\qquad\qquad\leq
\big(1+t \bar{d^2}/(2\underbar{d})\big)M_0.
\end{array}
\]
Finally, we estimate the negative part of the $a_i\ln(a_i)$:

\[\begin{array}{lll}
\ds\int_{\R^N} a_i|\ln(a_i)|\,dx
&=&
\ds\int_{\R^N} a_i\ln(a_i)\,dx
-2\ds\int_{\R^N} a_i\ln(a_i)\big(
1\!\!1_{0\leq a_i\leq e^{-|x|/2}}+
1\!\!1_{e^{-|x|/2}\leq a_i\leq 1}\big)\,dx
\\
&\leq&\ds\int_{\R^N} a_i\ln(a_i)\,dx
+ \ds\frac4e\ds\int_{\R^N} e^{-|x|/4} \,dx+
\ds\int_{\R^N} |x|a_i
\,dx
\end{array}\]
since $-s\ln(s)\leq \frac2e \sqrt s$ for any $0\leq s\leq 1$.
Combining together all the pieces we conclude the statement of the proposition.\\

\endproof

\vskip1cm

\section{Local parabolic regularization principle}\label{sec:bound}

We consider solutions of (\ref{eq_system}) that are defined for negative times. The following proposition is the main result of this section.

\vskip0.1cm

\begin{prop}\label{lemm_de giorgi}
Let (\ref{eq_system}) be a system verifying (\ref{Hypothesis_positivity}), (\ref{Hypothesis_subquadratic}), (\ref{Hypothesis_mass}), and (\ref{Hypothesis_entropy}). 
Let $a=(a_1,\cdot\cdot\cdot, a_P)$ be a solution to (\ref{eq_system}) on $(-3,0)\times B(0,3)$, such that $a_i$ are nonnegative functions for $1\leq i\leq P$. Then,
for any $p>1$, there exists a universal constant $\delta^{\star}>0$ (depending on
the constant $\Lambda$ as in (\ref{Hypothesis_subquadratic}), $\bar d$, $\underbar d$, $P$, $N$ and $p$), such that, if  $a=(a_1,\cdot\cdot\cdot, a_P)$
verifies
$$ \ds\sum_{i=1}^P \|a_i\|_{L^{ p}((-3,0)\times B(0,3))}\leq \delta^{\star},$$
then, $0\leq a_i(0,0)\leq 1$ for $\,1\leq i\leq P$.
\end{prop}
\vskip0.1cm
 The value at $(0,0)$ can be hence controlled by any $L^p$ norm, $p>1$, on a surrounding cylinder. Such results is quite surprising, since we work with subquadratic reaction terms.\\
\vskip0.1cm
In the spirit of the Stampacchia cut-off method, $L^\infty$ bounds
of solutions of certain PDEs can be deduced from the behavior of suitable non linear functionals.
Here, such functionals are constructed in a way that they use the dissipation property~\eqref{Hypothesis_entropy}. This is based on the De Giorgi techniques and it is reminiscent of the method introduced by Alikakos \cite{Alikakos}.

\vskip1cm
Let  us consider the non negative, $C^1$,  and convex function
\[\Phi(z)=\left\{\begin{array}{ll}(1+z)\ln(1+z)-z\qquad &\text{ if } z\geq 0,\\
0 \qquad&\text{ if } z\leq 0.\end{array}\right.\]
We study the evolution of the ''entropy at level $R$"
\[\ds\sum_{i=1}^P\ds\int_{ \mathcal{B}_n} \Phi(a_i-R)\,dx, \]
 for $R\geq 0$, where the set $ \mathcal{B}_n$ will be specified later.
\vskip1cm

\begin{lemm}\label{rineq} Let $1<\nu<2$. For any $a\in \R^P$, and for any $0\leq R\leq 1$,  we have
$$\sum_{i=1}^P\left | Q_i(a)-Q_i(1+[a-R]_+)\right |\leq 2 P\Lambda |1+[a-R]_+|^{\nu-1}$$
\noindent where $\Lambda$  is given by~(\ref{Hypothesis_subquadratic}) and $1+[a-R]_+:=(1+[a_1-R]_+, \cdot\cdot\cdot, 1+[a_P-R]_+)$.
\end{lemm}

\begproof
\noindent For $|u|\leq |v|$, by the sub-quadratic growth for the reaction term $Q$, hypothesis~(\ref{Hypothesis_subquadratic}), we have
$$\sum_{i=1}^P\left | Q_i(u)-Q_i(v)\right |\leq  \Lambda P |v-u| |v|^{\nu-1}.$$
\noindent The inequality follows by substituting $(u,v)$ by $(a,1+[a-R]_+)$.
\endproof
\vskip1cm
\noindent Let us set $k_n=1-1/2^n$, $t_n=1+1/2^n$, $\mathcal{B}_n=B(0,t_n), \mathcal{Q}_n=(-t_n,0)\times \mathcal{B}_n$.
\noindent Note that  $ \mathcal{B}_n\subset  \mathcal{B}_{n-1}$ and $ \mathcal{Q}_n\subset  \mathcal{Q}_{n-1}$. \\

We introduce the cut-off functions
\[\left\{\begin{array}{ll}
\zeta_n:\mathbb R^N\rightarrow \mathbb R,\qquad&
0\leq \zeta_n(x)\leq 1,\\
\zeta_n(x)=1 \textrm{ for $x\in \mathcal B_n$},\qquad&
 \zeta_n(x)=0 \textrm{ for $x\in\complement \mathcal B_{n-1}$},\\
 \ds\sup_{i,j\in \{1,\dots, N\},\ x\in\mathbb R^N} | \partial_{ij}^2\zeta_n(x)|\leq C\  2^{2n}\,\,\,\hbox{with $C$ universal constant}. \qquad&
\end{array}\right.\]

Next, we give an estimate on the local dissipation of entropy at the level $R$.

\begin{prop}\label{gainofr}
There exists a universal constant $\hat C$ (depending on $\bar d$,$\underbar d$, $\Lambda$ and $P$), such that for every $a=(a_1,\cdot\cdot\cdot, a_P)$ solution of~(\ref{eq_system}), for any $0\leq R\leq 1$, we have
$$\sup_{-t_n\leq t\leq0}\big\{\sum_{i=1}^{P}\,\int_{\mathcal B_n} \Phi(a_i-R)(t,x)\,dx\big\}+ \underbar{d}\sum_{i=1}^{P}\,\int \int_{\mathcal Q_n}|\nabla_x\sqrt{1+[a_i-R]_+}|^2\,dx\,d\tau$$
$$\leq 2^{2n}\hat C \sum_{i=1}^{P}\int\int_{\mathcal Q_{n-1}}(1+[a_i-R]_+)\ln(1+[a_i-R]_+)\,dx\,d\tau.$$
\end{prop}

\begproof
We multiply (\ref{eq_system}) by $\zeta_n\Phi'(a_i-R)$, and we sum
\begin{equation}\label{int1}
\ds\frac{d}{dt}
\ds\sum_{i=1}^{P}\int_{\mathcal B_{n-1}}
 \zeta_n\Phi(a_i-R)\,dx= A+B
 \end{equation}
 \noindent where
$$ A:= \sum_{i=1}^{P} \int_{\mathcal B_{n-1}}
 D_i\Delta a_i\
 \Phi'(a_i-R)\zeta_n\,dx$$
 \noindent and
 $$B:= \ds\sum_{i=1}^{P}\int_{\mathcal B_{n-1}}
 Q_i(a)\Phi'(a_i-R)\zeta_n\,dx.$$

 We rewrite $A$ as,
 $$
A=  -E+F
 $$
  \noindent where

 $$
 E:= \ds \sum_{i=1}^{P} \int_{\mathcal B_{n-1}}
 D_i|\nabla a_i|^2\
 \Phi''(a_i-R)\zeta_n\,dx$$
  \noindent and
 $$F:= \sum_{i=1}^{P} \int_{\mathcal B_{n-1}}
 D_i \Phi(a_i-R)\Delta\zeta_n\,dx
. $$
 Moreover,
 \[\begin{array}{l}
E= \ds\sum_{i=1}^{P}\int_{\mathcal B_{n-1}}
 D_i|\nabla_x a_i|^2\
 \Phi''(a_i-R) \zeta_n\,dx=
 \ds\sum_{i=1}^{P}\int_{\mathcal {B}_{n-1}}
 D_i\nabla_x a_i\cdot\nabla_x a_i\
 \ds\frac{1\!\!1_{a_i\geq k}}{1+[a_i-R]_+}\,dx
 \\
 \qquad=
 \ds\sum_{i=1}^{P}\int_{\mathcal B_{n-1}}
 D_i\nabla_x (1+[a_i-R]_+)\cdot\nabla_x (1+[a_i-R]_+)\
 \ds\frac{\,dx}{1+[a_i-R]_+}
 \\
\qquad \geq
 \underline{d}\ds\sum_{i=1}^{P}\int_{\mathcal B_{n-1}}
 \ds\frac{\big|\nabla_x (1+[a_i-R]_+)\big|^2}{1+[a_i-R]_+}
 \,dx
 \\
\qquad \geq
 4\underline{d}
 \ds\sum_{i=1}^{P}\int_{\mathcal B_{n-1}}\big|\nabla_x \sqrt{1+[a_i-R]_+}\big|^2\,dx.
 \end{array}\]
Next, we estimate the quantity $B$ from~\eqref{int1}. For $0<\nu\leq 1$,
$$
 B\leq \ds\sum_{i=1}^{P}
 \ds\int_{\mathcal B_{n-1}}
 (1+[a_i-R]_+) \ln(1+[a_i-R]_+)\,dx,
 $$
 while, for $1<\nu<2$, to get rid of the nonlinearity given by the factor $Q_i(a)$, we
rewrite $B$ as 
 \[
 \begin{array}{l}
 \ds\sum_{i=1}^{P}
 \ds\int_{\mathcal B_{n-1}}
 Q_i(a)\
 \ln(1+[a_i-R]_+)\,dx
 \\
 =
  \ds\sum_{i=1}^{P}
 \ds\int_{\mathcal B_{n-1}}
 \big(Q_i(a)-Q_i(1+[a-R]_+)\big)\
 \ln(1+[a_i-R]_+)\,dx
 \\
 \qquad\qquad+
 \ds\sum_{i=1}^{P}
 \ds\int_{\mathcal B_{n-1}}
 Q_i(1+[a-R]_+)
 \ln(1+[a_i-R]_+)\,dx.
 \end{array}
 \]
By assumption~\eqref{Hypothesis_entropy} the last term is non positive,
while $ \big(Q_i(a)-Q_i(1+[a-R]_+)\big)\
 \ln(1+[a_i-R]_+)$ can be estimated via Lemma~\ref{rineq}. Finally, from~(\ref{int1}) and the above, we get
\begin{eqnarray*}
&&\ds\frac{d}{dt}
\ds\sum_{i=1}^{P}\int_{\mathcal B_{n-1}}
 \zeta_n\Phi(a_i-R)\,dx
 +
4\underline{d}
 \ds\sum_{i=1}^{P}\int_{\mathcal B_n}\big|\nabla_x \sqrt{1+[a_i-R]_+}\big|^2\,dx\\
&&\qquad
\leq \hat C\left(\ds\sum_{i=1}^{P}\int_{\mathcal B_{n-1}}
(1+[a_i-R]_+)\ln (1+[a_i-R]_+)\,dx+ \ds \sum_{i=1}^{P}\int_{\mathcal B_{n-1}} \Phi(a_i-R)\,dx\right).
\end{eqnarray*}

To conclude, we integrate on $(s,t)$, for $-t_n<t<0$, and average on $s\in(-t_{n-1},-t_n)$.
\endproof
\vskip0.1cm

Let us introduce the following easy lemma.
\begin{lemm}\label{easy}
Let
\begin{equation}\label{defPsi}
\Psi(z)=(\sqrt{1+z}-1)\un_{\{z>0\}},\qquad z\in \R.
\end{equation}
The function $\Psi$ is nondecreasing, Lipschitz  on $\R$, and there exists a constant $C>0$ such that
\begin{eqnarray}
&&\Psi(z)\leq \tilde C\sqrt{\Phi(z)},\qquad z\in \R,\label{Psi controled by phi}\\
&&\un_{\{\Psi(z-k_n)>0\}}\leq \tilde C 2^n \psi(z-k_{n-1}),\qquad n\geq 1.\label{De Giorgi trick on psi}
\end{eqnarray}
\end{lemm}

\noindent
\begproof
To get (\ref{Psi controled by phi}), we expand  $\Psi$ and $\Phi$ both at $0$ and $+\infty$ to find that there exists $C_1, C_2>0$ with
\begin{eqnarray*}
\Psi(z)^2\leq C_1 \inf(z,z^2)\qquad z>0,\\
\Phi(z)\geq C_2 \inf(z,z^2)\qquad z>0.
\end{eqnarray*}
To get (\ref{De Giorgi trick on psi}), first consider $z>2$:\\
For any $n\geq 1$
$$
2^n\Psi(z-k_{n-1})\geq \Psi(z-k_{n-1})\geq \Psi(1)\geq (\sqrt{2}-1)\un_{\{\Psi(z-k_n)>0\}}.
$$
Next, if $z\leq 2$, $\Psi(z-k_n)>0$ implies that $z>k_n$ and so $z-k_{n-1}\geq 2^{-n}$.
Hence, for such $z$:
\begin{eqnarray*}
\Psi(z-k_{n-1})&\geq&[\Psi(z-k_{n-1})-\Psi(z-k_{n})]+\Psi(z-k_{n})\\
&\geq&\Psi(z-k_{n-1})-\Psi(z-k_{n})\\
&\geq&\frac{k_n-k_{n-1}}{\sqrt{1+(z-k_{n-1})_+}+\sqrt{1+(z-k_{n})_+}}\geq \frac{2^{-n}}{2\sqrt{3}}.
\end{eqnarray*}

\endproof

\noindent We, now, define the sequence $\mathcal{U}_n$ which plays a key role for the proof of the uniform boundedness of the solution.\\

$$\mathcal{U}_n=\sup_{-t_n\leq t\leq 0}\sum_{i=1}^{P}\int_{{\mathcal{B}}_n} \Phi(a_i-k_n)\,dx+ \sum_{i=1}^{P}\int\int_{{\mathcal{Q}}_n}|\nabla_x\Psi(a_i-k_n)|^2\,dx\,ds.$$
Note that, since $\Psi$ is Lipschitz,
$$
\nabla_x\Psi(a_i-k_n)=\nabla_x\sqrt{1+[a_i-k_n]_+}.
$$

The following nonlinear estimate is a crucial step for establishing our results.\\

\begin{lemm}\label{unestimates} There exists a universal constant $C$  (depending only on $\Lambda$, $\bar d$, $\underbar d$, $P$, and $N$)
such that
$$\mathcal{U}_n\leq C^n\mathcal{U}_{n-1}^{\frac{N+2}{N}}$$
\noindent for any $n\geq 1$. Especially, there exists $\delta>0$   (depending only on $\Lambda$, $\bar d$, $\underbar d$, $P$, and $N$) such that if  $\mathcal{U}_0\leq \delta$ then $\lim_{n->\infty} \mathcal{U}_n=0$.

\end{lemm}

\noindent
\begproof
\noindent From Proposition~\ref{gainofr} we get the
 inequality:
\begin{eqnarray*}
&&\mathcal U_n\leq \hat C 2^{2n}\,
\
\ds\sum_{i=1}^{P}\ds\int\int_{\mathcal{Q}_{n-1}}
\ln(1+[a_i-k_{n}]_+)(1+[a_i-k_n]_+) \,dx\,d\tau\\
&&\qquad
\leq \hat C2^{2n}\ds\sum_{i=1}^{P}\ds\int\int_{\mathcal{Q}_{n-1}}\un_{\{a_i\geq k_n\}}
(1+[a_i-k_n]_+)^{\frac{N+2}{N}} \,dx\,d\tau.
\end{eqnarray*}
Hence
\begin{equation}\label{crucial2}
\mathcal U_n\leq 2\hat C 2^{2n}\left(\sum_{i=1}^{P}\ds\int\int_{\mathcal{Q}_{n-1}}\un_{\{\Psi(a_i-k_n)>0\}}\,dx\,d\tau+\|\Psi(a_i-k_n)\|^{2\frac{N+2}{N}}_{L^{2\frac{N+2}{N}}(Q_{n-1})}\right).
\end{equation}

Next, using (\ref{De Giorgi trick on psi}), we find that
\begin{eqnarray*}
&&\qquad\qquad \sum_{i=1}^{P}\ds\int\int_{\mathcal{Q}_{n-1}}\un_{\{\Psi(a_i-k_n)>0\}}\,dx\,d\tau\\
&&\leq (\tilde C2^n)^{2\frac{N+2}{N}}\sum_{i=1}^{P}\ds\int\int_{\mathcal{Q}_{n-1}}[\Psi(a_i-k_{n-1})]^{2\frac{N+2}{N}}\,dx\,d\tau.
\end{eqnarray*}
Since $\Psi(a_i-k_n)\leq \Psi(a_i-k_{n-1})$, we finally get

$$
\mathcal U_n\leq  2\hat C 2^{2n} (\tilde C2^n)^{2\frac{N+2}{N}}\sum_{i=1}^{P}\,\|\Psi(a_i-k_{n-1})\|^{2\frac{N+2}{N}}_{L^{2\frac{N+2}{N}}(\mathcal{Q}_{n-1})}.
$$
Using (\ref{Psi controled by phi}), we find that
$$
\|\Psi(a_i-k_{n-1})\|^2_{L^\infty(-t_{n-1},0;L^2(\mathcal{B}_{n-1}))}\leq \tilde C \|\Phi(a_i-k_{n-1})\|_{L^\infty(-t_{n-1},0;L^1(\mathcal{B}_{n-1}))}\leq \mathcal U_{n-1}.
$$
By Sobolev imbedding, we find
$$
\|\Psi(a_i-k_{n-1})\|^2_{L^2(-t_{n-1},0;L^{\frac{2N}{N-2}}(\mathcal{B}_{n-1}))}\leq
c \,\mathcal U_{n-1}\,\,\hbox {}.
$$
By interpolation, this gives
$$
\|\Psi(a_i-k_{n-1})\|^2_{L^\frac{2(N+2)}{N}(\mathcal{Q}_{n-1})}\leq
c \mathcal U_{n-1} ,\,\hbox {with $c$ universal constant}.
$$
Finally, we get
$$
\mathcal U_n\leq C^n {\mathcal U}^\frac{(N+2)}{N}_{n-1}
$$
where $C$ depends only on
the constants $\Lambda$ as in~(\ref{Hypothesis_subquadratic}), $\bar d$, $\underbar d$, $P$, and $N$.
\endproof

We now discuss the smallness of~$\mathcal U_0$.

\begin{lemm}\label{lem U_0} For any  $ p>1$ there exists  a universal constant $C>0$ (depending only on $\bar d$, $\underbar d$ and $p$) such that
$$\mathcal U_0\leq C\big ( \ds\sum_{i=1}^{P} \|a_i\|^p_{L^{ p}((-3,0)\times  B(0,3))}+  \ds\sum_{i=1}^{P} \|a_i\|^{1/2}_{L^{ p}((-3,0)\times B(0,3))}\big ).$$

 \end{lemm}
\begproof
\noindent By definition
\[\mathcal U_0=
\ds\sup_{-2\leq t\leq 0}
\ds\sum_{i=1}^{P}\ds\int_{ B(0,2)}
\Phi(a_i)\,dx
+
\ds\int_{-2}^0\ds\int_{B(0,2)}\big|\nabla_x\sqrt{a_i+1}\big|^2\,dx\,d\tau.
\]

\noindent We will use the following facts:\\
\begin{align}
& \Phi(a_i)\leq C (a_i (1+|\ln(a_i)|)\label{fact2};\\
& \ds\int_{-2}^0\ds\int_{B(0,2)}\big|\nabla_x\sqrt{a_i+1}\big|^2\,dx\,d\tau\leq \ds\int_{-2}^0\ds\int_{B(0,2)}\big|\nabla_x\sqrt{a_i}\big|^2\,dx\,d\tau. \label{fact1}
\end{align}
\noindent The mass conservation yields
\[\ds\frac{d}{dt}\ds\sum_{i=1}^{P}\ds\int_{B(0,3)} \zeta_0(x)a_i(t,x)\,dx
=
\ds\sum_{i=1}^{P}\ds\int_{B(0,3)} D_i \Delta\zeta_0(x)\ a_i(t,x)\,dx
\leq C \ds\sum_{i=1}^{P}\ds\int_{B(0,3)} a_i(t,x)\,dx,\]
where $\zeta_0\geq 0$,    $\zeta_0\equiv 1$ on $B(0,2)$,  has bounded second order derivatives and it is supported on $B(0,3)$.
Let $t\in(-2, 0)$. Let $\tau\in (-3,t)$. We integrate over the time interval $(\tau,t)$, and then we average over $\tau\in (-3,-2)$.
Hence, we get
\[\begin{array}{lll}
\ds\sup_{-2\leq t\leq 0}\ds\sum_{i=1}^{P}\ds\int_{B(0,2)}a_i(t,x)\,dx
&\leq&
C\ds\sum_{i=1}^{P} \ds\int_{-3}^0\ds\int_{B(0,3)}
a_i(\tau,x)\,dx \,d\tau.

\end{array}\]\\
Similarly, the entropy dissipation yields
\[\begin{array}{l}
\ds\frac{d}{dt}\ds\sum_{i=1}^{P}\ds\int_{B(0,3)} \zeta_0(x)\ a_i\ln (a_i)\,dx
+ \ds\int_{B(0,3)} \zeta_0(x)\ds\frac{D_i |\nabla_x{a_i}|^2}{a_i}\,dx
\\
\qquad\leq
\ds\sum_{i=1}^{P}\ds\int_{B(0,3)} D_i \Delta\zeta_0(x)\ a_i\ln (a_i)\,dx
\\
\qquad
\leq C \ds\sum_{i=1}^{P}\ds\int_{B(0,3)} a_i|\ln(a_i)|\,dx.
\end{array}\]
Again we  integrate with respect to the time variable.
We shall also use the trick
\[
u|\ln(u)|=
u\ln(u)-2u\ln(u)1\!\!1_{0\leq u\leq 1}\leq
u\ln(u)+C\sqrt{u}.
\]
It follows that
\[\begin{array}{l}
\ds\sup_{-2\leq t\leq 0}\ds\sum_{i=1}^{P}\ds\int_{B(0,2)}a_i|\ln(a_i)|\,dx
+\underline{d} \ds\sum_{i=1}^{P} \ds\int_{-2}^0\ds\int_{B(0,2)}\big|\nabla_x\sqrt{a_i}\big|^2\,dx\,d\tau
\\
\qquad\qquad\leq
 C \ds\sum_{i=1}^{P}\left ( \ds\int_{-3}^0\ds\int_{B(0,3)}
\left(a_i|\ln(a_i)|+\sqrt{a_i}\right)\,dx\,d\tau\right).
\end{array}\]


\noindent Combining~(\ref{fact1}),  ~(\ref{fact2}), and the definition of $\mathcal U_0$ yields to the desired inequality.
\endproof
\vskip1cm
The proof of Proposition \ref{lemm_de giorgi} follows from Lemma \ref{unestimates} and Lemma \ref{lem U_0}.

\section{Duality arguments}\label{sec:duality}

In this section we derive a uniform bound on  the solutions $a=(a_1,\cdot\cdot\cdot, a_P)$ in a weak norm. We will show that this weak norm shrinks through the universal scaling. The proof relies on a nice duality argument first used in \cite{PiSc}.
 
 Let $a=(a_1,\cdot\cdot\cdot, a_P)$ be a solution to the system~(\ref{eq_system}) with initial data $a^0=(a_1^0,\cdot\cdot\cdot, a_P^0)$. Then, the function $\rho=\sum_{i=1}^P a_i$ is a solution to the Cauchy problem
\begin{equation}\label{eq_laplace}
\left\{\begin{array}{l}
\displaystyle{\dt \rho-\Delta(d\rho)=0, \qquad t>0, \ x\in
\R^N,}\\[0.3cm]
\displaystyle{\rho^0(x)=\rho(0,x)=\sum_{i=1}^Pa^0_i,
\qquad x\in \R^N,}
\end{array}
\right.
\end{equation}
with
$$d(t,x)=\frac{\sum_{i=1}^P D_i a_i}{\sum_{i=1}^Pa_i}.$$

 Note that the diffusion coefficient $d(t,x)$  is elliptic. Indeed, it is bounded from above and below:
\begin{equation}\label{borne sur d}
0<\underbar{d}\leq d(t,x)\leq \bar{d}.
\end{equation}
The equation in~(\ref{eq_laplace}) seems so to be a nice parabolic equation, except that it is not in the standard divergence or non divergence form. In the non divergence form, the equation would provide the maximum principle. In the classical divergence form, De Giorgi showed in \cite{DeG} that such solutions are bounded  locally in $C^\alpha$.
Surprisingly, the behavior of solutions of parabolic equations written as (\ref{eq_laplace}) is very different. Note that we do not have a priori bounds based on the regularity of $d$. In \cite{PiSc}, Pierre and Schmitt show that any solution of the parabolic equation~(\ref{eq_laplace}), with (\ref{borne sur d}) (but no assumption on the regularity of $d$), and  with regular enough initial values, lies in $L^2((0,T_0)\times\R^N)$ for $T_0>0$. However, they also give explicit examples of solutions which blow-up in $L^p((0,T_0)\times\R^N)$ for some $p>2$. It follows that global regularity of  solutions to (\ref{eq_system}) cannot
 rely only on the equation solved by $\rho$.

As usual we denote $ \mathcal{D}(\R^N):=\mathcal{C}^{\infty}_0(\R^N)$. Let us first introduce our weak space.
\begin{defi} We define the space $\mathcal{L}_w(\R^N)$ as the dual space of $\{\rho\in \mathcal{D}(\R^N): \Delta\rho\in L^1\}$:

$$\mathcal{L}_w(\R^N):=\{f\in \mathcal{D}^{\prime}(\R^N): ||f||_{\mathcal{L}_w(\R^N)}<\infty\}$$
\noindent where
$$ ||f||_{\mathcal{L}_w(\R^N)}:=\sup_{\psi \in \mathcal{D}(\R^N),
\left||\Delta\psi \right \|_{L^1}\leq 1}\left |<f,\psi>\right |.$$
\end{defi}


The following proposition is the main result of this section.
\begin{prop}\label{lemm_W2infini}
Let $a=(a_1,\cdot\cdot\cdot, a_P)$ be a smooth and bounded solution of (\ref{eq_system}) on the time interval $[0,T]$ for any $T<T_0$ (with possible blow-up at $T=T_0$). Then $\rho=\sum_{i=1}^P a_i$ verifies $$\rho\in L^\infty(0,T_0;\mathcal{L}_w(\R^N)),\,\,\,$$
\noindent with,
\begin{equation}
\|\rho(T)\|_{\mathcal{L}_w(\R^N)}\leq \|\rho^0\|_{\mathcal{L}_w(\R^N)},\,\,\,\,\hbox{for all} \,\,T<T_0.\label{uniboundweak}
\end{equation}
\end{prop}
\begproof Let us fix $T<T_0$. The proof consists of two steps:

\vskip0.1cm
\noindent  {\it Step 1}: Since $a=(a_1,\cdot\cdot\cdot, a_P)$ is smooth, the function $d$ is smooth also, except possibly at the points $(t,x)$ where $\rho=0$. For small $\mu>0$, let $d_{\mu}$ denote a smooth approximation of $d$ verifying
\begin{equation}\label{prop d}
\|d_\mu\|_{L^\infty((0,T)\times\R^N)}\leq \overline{d},\qquad d_\mu(t,x)=d(t,x)\ \ \mathrm{if} \ \ \rho(t,x)\geq \mu.
\end{equation}
\noindent Let $\rho_{\mu}$ be a solution to the problem:
\begin{equation}\label{mueq_laplace}
\begin{array}{l}
\displaystyle{\dt \rho_{\mu}-\Delta(d_{\mu}\rho_{\mu})=0,\,\,\, x\in\R^N,\,\,\, 0<t<T,}\\
\displaystyle{\rho(0,\cdot)=\rho_{\mu}(0)\in \mathcal{D}(\R^N)\cap L^1(\R^N).}
\end{array}
\end{equation}
\vskip0.2cm
We claim that for such a solution, we have for every $0\leq t\leq T$:
\begin{equation}
\|\rho_{\mu}(t)\|_{L^1(\R^N)}\leq \|\rho_{\mu}(0)\|_{L^1(\R^N)}.\label{unil1bound}
\end{equation}
To show this claim, consider the dual problem, for any $0<{\overline{T}}<T$: 
\begin{equation}\label{eq_dual1}
\begin{array}{l}
\displaystyle{\dt\phi+d_{\mu}\Delta\phi=0,\qquad  x\in \R^N, 0\leq t\leq {\overline{T}},}\\
\displaystyle{\phi({\overline{T}})=\phi_{\overline{T}}\in \mathcal{D}(\R^N)\cap L^\infty(\R^N).}
\end{array}
\end{equation}
\noindent The maximum principle, applied to the problem~(\ref{eq_dual1}), gives a uniform bound for the $L^{\infty}$ norm of $\phi$:
$$\|\phi(t)\|_{L^{\infty}(\R^N)}\leq \|\phi_{\overline{T}}\|_{L^{\infty}(\R^N)},$$
for any $0<t<{\overline{T}}$.
But, for every $0<t<{\overline{T}},$
$$
\frac{d}{dt}\int_{\R^N}\rho_\mu(t,x)\phi(t,x)\,dx=0.
$$
Then, for any $\phi_{\overline{T}}\in\mathcal{D}(\R^N)$ with $\|\phi_{\overline{T}}\|_{L^\infty(\R^N)}\leq 1$, we have $\|\phi(0)\|_{L^\infty(\R^N)}\leq 1$ and
\begin{eqnarray*}
\left|\int_{\R^N}\rho_\mu({\overline{T}},x)\phi_{\overline{T}}(x)\,dx\right|&=&  \left|\int_{\R^N}\rho_\mu(0,x)\phi(0,x)\,dx\right|\\
&\leq&\|\rho_\mu(0)\|_{L^1(\R^N)}\|\phi(0)\|_{L^\infty(\R^N)}\\
&\leq&\|\rho_\mu(0)\|_{L^1(\R^N)}.
\end{eqnarray*}
This proves claim (\ref{unil1bound}).

\vskip0.1cm
\noindent {\it Step 2}: We consider, now, solutions to the following dual problem 
\begin{equation}\label{eq_dual2}
\begin{array}{l}
\displaystyle{\dt\phi+d_{\mu}\Delta\phi=0,\qquad  x\in \R^N, \,\,0\leq t\leq T,}\\
\displaystyle{\phi(T)=\phi_{T}\in \mathcal{D}(\R^N),}
\end{array}
\end{equation}
for any $\phi_T$ such that 
$$
\|\Delta \phi_T\|_{L^1(\R^N)}\leq 1.
$$
The function
$$
\rho_\phi(t,x)=\Delta \phi(T-t,x), \qquad 0\leq t\leq T, x\in \R^N,
$$
verifies (\ref{mueq_laplace}) with initial value $\Delta \phi_T$.
Hence, (\ref{unil1bound}) ensures that
\begin{equation}\label{laplace}
\|\Delta\phi(t)\|_{L^1(\R^N)}\leq \|\Delta\phi_T\|_{L^1(\R^N)}\leq 1\qquad
0\leq t\leq T.
\end{equation}
But, for any solution $a=(a_1,\cdot\cdot\cdot, a_P)$ of (\ref{eq_system}),  $\rho=\sum_{i=1}^P a_i$ verifies (\ref{eq_laplace}) and so we have:
$$\frac{d}{dt}\int_{\R^N} \rho(t,x)\phi\,dx=\int_{\R^N} \rho(t,x)(d-d_{\mu})\Delta\phi\,dx.$$
Integrating in time on $[0,T]$, we find
\begin{eqnarray*}
&&\left|\int_{\R^N}\phi_T(x)\rho(T,x)\,dx\right|\leq\left|\int_{\R^N}\phi(0,x)\rho^0(x)\,dx\right|+\left|\int \int_{\R^N} \rho(d-d_{\mu})\Delta\phi\,dx\,dt\right|\\
&&\qquad\leq \|\rho^0\|_{\mathcal{L}_w(\R^N)}+2 T\mu\overline{d}\|\Delta\phi\|_{L^{\infty}(0,T;L^1(\R^N))}\\
&&\qquad\leq \|\rho^0\|_{\mathcal{L}_w(\R^N)}+2 T\mu\overline{d}.
\end{eqnarray*}
In the second line we have used the definition of the weak norm together with (\ref{laplace}), and (\ref{prop d}). In the last line we have used again (\ref{laplace}). 

Passing to the limit, as $\mu\to0$,  for any $\phi_T\in \mathcal{D}(\R^N)$ such that $\|\Delta\phi_T\|_{L^1(\R^N)}\leq 1$, we have
$$
\left|\int_{\R^N}\phi_T(x)\rho(T,x)\,dx\right|\leq \|\rho^0\|_{\mathcal{L}_w(\R^N)}.
$$ 
This implies, thanks to the definition of the weak norm that:
$$
\|\rho(T)\|_{\mathcal{L}_w(\R^N)}\leq \|\rho^0\|_{\mathcal{L}_w(\R^N)}.
$$
The result follows, since this holds for any $T<T_0$.

\endproof
\vskip0.2cm
As a consequence of this proposition we have the following uniform estimate on the $\mathcal{L}_w(\R^N)$ norm of the total mass.

\begin{corr}\label{cor_ess}
There exists a universal constant $C$, such that, for any
 $a=(a_1,\cdot\cdot\cdot, a_P)$  regular and bounded solution to (\ref{eq_system}) on $[0,T]\times\R^N$ for any $T<T_0$ (with possible blow-up at $T_0$), we have $\rho\in L^\infty(0,T_0;\mathcal{L}_w(\R^N))$, and for any $t<T_0$
$$
\|\rho(t)\|_{\mathcal{L}_w(\R^N)}\leq C\left(\|\rho^0\|_{L^\infty({\R^N})}+\|\rho^0\|_{L^1{(\R^N)}}\right).
$$
\end{corr}

\begproof
In the view of Proposition \ref{lemm_W2infini}, we only need to show that the weak norm can be controlled by stronger norms. We set
$$
\phi=\Gamma\ast \rho^0
$$
where $\Gamma(x)=\frac{1}{N(N-2)\omega_N}\frac{1}{|x|^{N-2}}$ ($\omega_N$ is the volume of the unit ball in $\R^N$);
then, the function $\phi$ solves
$$
-\Delta\phi=\rho^0.
$$
Since $\Gamma\in L^1({\R^N})+L^\infty({\R^N})$, there exists a constant $C$ such that 

$$\|\phi\|_{L^{\infty}(\R^N)}\leq C(\|\rho^0|\|_{L^1({\R^N})}+\|\rho^0\|_{L^{\infty}({\R^N})}).$$
Therefore, for any $\psi\in\mathcal{D}(\R^N)$, with $\|\Delta\psi\|_{L^1(\R^N)}\leq 1$, we have
$$|\int_{\R^N} \rho^0\psi |=|\int_{\R^N} \Delta\phi\psi|=|\int_{\R^N} \phi\Delta\psi|\leq \|\phi\|_{L^{\infty}({\R^N})}\|\Delta\psi\|_{L^1({\R^N})}\leq C(\|\rho^0|\|_{L^1({\R^N})}+\|\rho^0\|_{L^{\infty}({\R^N})}).$$
\endproof

We observe that the nonnegativity of the mass $\rho$ allows us to control stronger norms from this weak norm.

\begin{lemm}\label{l1bound} Let $f$ be a nonnegative function on $\R^N$ such that $f\in L^{\infty}(\R^N)\cap L^1(\R^N)$, then
for any compact subset $K$ of $\R^N$,
$K\subset\subset \R^N$, there is a positive constant $C(K)$ such that
\begin{equation}\label{boundl1loc}
 \|f\|_{L^1(K)}\leq C(K) \|f\|_{\mathcal{L}_w(\R^N)}.
 \end{equation}
\end{lemm}
\begproof
Let $K$ be a compact set in $\R^N$ such that $K\subset \mathcal B(0,r)$, where $\mathcal B(x_0,r)$ is a ball centered at some $x_0\in K$ of radius $r>0$ . Let $\phi\in C^\infty(\R^N)$ be a smooth nonnegative function, such that $\phi = 1$ on $K$,  and  $\phi=0$ on  $\complement  (\mathcal B(x_0,2r))$. Set
$$
C(K)=\|\Delta\phi\|_{L^1(\R^N)}.
$$
Then, the $L^1$ norm of the Laplacian of $\displaystyle{\frac{\phi}{C(K)}}$  is equal to 1, which implies that
$$
\|f\|_{L^1(K)}\leq \int_{\R^N}\phi(x)f(x)\,dx\leq C(K)\int_{\R^N}\frac{\phi(x)}{C(K)}f(x)\,dx\leq C(K)\|f\|_{\mathcal{L}_w(\R^N)}.
$$
\endproof

As a consequence of this lemma we have the following corollary.

\begin{corr}\label{corrlqlp} For any integer $P>0$, and any $p\in(1,\frac{N}{N-2})$, there exists a universal constant $C$ (depending only on $p$, $P$ and $N$) such that the following is true. Given nonnegative functions $a_i=a_i(t,x)$, for $1\leq i\leq P$,  such that
$$
\rho=\sum_{i=1}^Pa_i\in L^\infty(-3,0;\mathcal{L}_w(\R^N)),\qquad  \ \ \nabla\sqrt{a_i}\in L^2((-3,0)\times\R^N),\,\,\, 1\leq i\leq P,
$$
we have
$$
 a_i\in L^q(-3,0;L^p(B(0,3))),\,\,\,\,\,\,\,\,\,\,
 $$
 \noindent  for $q$ such that  $\frac{1}{p}=1-\frac{2}{qN}$.
 Moreover,
 \begin{equation}\label{mainpropertyforscaling}
\|a_i\|_{L^q(-3,0;L^p(B(0,3)))}\leq C \|\nabla_x\sqrt  a_i \|_{ L^2((-3,0)\times R^N)}^{\frac{2(p-1)}{p}}\|\rho\|_{L^\infty(-3,0;\mathcal{L}_w(\R^N))}^{\frac{1}{p}}.
\end{equation}
\end{corr}
\begproof 
By Lemma~\ref{l1bound},  there exists a constant $C$ such that the inequality~(\ref{boundl1loc}) holds
$$ ||\rho(t,\cdot)||_{L^1(B(0,3))}\leq C ||\rho(t,\cdot)||_{\mathcal{L}_w(\R^N)}\,\,\,\,\hbox{for} \,\,-3<t<0.$$
And so, using the nonnegativity of the functions, for any $1\leq i\leq P$, we have
$$
\|a_i\|_{L^\infty(-3,0;L^1(B(0,3)))}\leq C ||\rho||_{L^\infty(-3,0;\mathcal{L}_w(\R^N))}.
$$
Moreover, by the Sobolev imbedding
$$
\|a_i\|_{ L^1(-3,0; L^{\frac{N}{N-2}} (R^N))}\leq C\|\nabla_x\sqrt  a_i \|_{ L^2((-3,0)\times R^N)}. 
$$
Now,  let $0<\theta<1$,  and $q$ defined by
$$
    \left\{
    \begin{aligned}\frac{1}{p}= & \frac{1-\theta}{1}\\
 \frac{1}{q}= & \frac{\theta}{1}+\frac{1-\theta}{\frac{N}{N-2}}
   \end{aligned}
\right.
$$
Then, by standard interpolation, it follows that  each $a_i$ belongs to the space $ L^q(0,T_0;L^p(B(0,3)))$, with the bound
\begin{equation}\label{boundinterpolation}
\|a_i\|_{L^q(-3, 0;L^p(B(0,3)))}\leq C(B(0,3)) \|a_i\|_{L^1(-3,0;L^{\frac{N}{N-2}}(\R^N))}^{\theta}
\|a_i\|_{L^{\infty}(-3,0;L^1(B(0,3)))}^{1-\theta}.
\end{equation}

\noindent The inequality~(\ref{mainpropertyforscaling}) follows readily by the definition of $\rho$.
\endproof

\section{Scaling argument}

This section is dedicated to the proof of Theorem \ref{Theo_principal}.\\
\vskip0.1cm
The proof is done by a contradiction argument. First, by standard short time existence results, we know that there exists a solution for some short time. Next, we assume that the maximal time of existence, which we denote by $T_0$, is finite.
We will show that the $L^{\infty}$ norm of $a=(a_1,\cdot\cdot\cdot, a_P)$ is uniformly bounded on $(\frac{T_0}{2},T_0)\times \R^N$. But this will contradict the fact that the solution blows up at $T_0$.\\
\vskip0.5cm
We introduce now the rescaled solutions. Let $a=(a_1,\cdot\cdot\cdot, a_P)$ be a solution of~(\ref{eq_system}).  Let $\frac{T_0}{2}<T<T_0$, $x_0\in \R^N$ and $0<\eps<\sqrt{T_0/6} $. Then, we define the rescaled functions
\begin{equation}\label{eq_a epsilon}
a^{\eps}(s,y)=\eps^{\frac{2}{\nu-1}}a(\eps^2 s+T,\eps y+x_0), \,\, -\frac{T}{\eps^2}\leq s\leq \frac{T_0-T}{\eps^2},\,\,\, y\in\R^N.
\end{equation}

\noindent The function $a^\eps=(a_1^\eps,\cdot\cdot\cdot, a_P^\eps)$  is a solution to a system (\ref{eq_system}) with
$$
Q^{\eps}(a)=\eps^{\frac{2\nu}{\nu-1}}Q\left(\eps^{\frac{-2}{\nu-1}}a\right).
$$
The reaction term $Q^{\eps}(a)$ verifies Hypothesis
(\ref{Hypothesis_positivity}), (\ref{Hypothesis_subquadratic}),
 (\ref{Hypothesis_mass}), and (\ref{Hypothesis_entropy})
 with
the constant $\Lambda$ as in (\ref{Hypothesis_subquadratic}) 
independent on $\eps$.
The $\mathcal{L}_w$ norm of the function $\rho^{\eps}$ rescales by the following identity
$$
\|\rho^{\eps}\|_{L^\infty(-\frac{T}{\eps^2}, 0;\mathcal{L}_w(\R^N))}=\eps^{\frac{2}{\nu-1}-2}\|\rho\|_{L^\infty(0,T;\mathcal{L}_w(\R^N))}.
$$
This implies,  by Corollary~\ref{cor_ess}, that the following bound holds
\begin{equation}\label{final control weak norm}
\|\rho^{\eps}\|_{L^\infty(-3, 0;\mathcal{L}_w(\R^N))}\leq C\eps^{\frac{2}{\nu-1}-2}(\|\rho^0\|_{L^1(\R^N)}+\|\rho^0\|_{L^\infty(\R^N)}).
\end{equation}
Through the rescaling, from~(\ref{final control weak norm}), we can control the following $L^{\bar q}$ norm. 

\begin{lemm}\label{lemm_de giorgi2} There exists $\bar q>1$ and $\eps_0>0$ (depending on $\bar d$, $\underbar d$, $\Lambda$, $P$, $N$, $T_0$ and $M_0$)  such that for all $\eps\leq\eps_0$ we have
$$\sum_{i=1}^{P} \|a^{\eps}_i\|_{L^{\bar q}((-3,0)\times B(0,3))}\leq \delta^{\star},$$
where $\delta^{\star}=\delta^{\star}(\bar q)$ is the same as in Proposition~\ref{lemm_de giorgi}.
\end{lemm}
\begproof
\noindent Let $0<\eps_0<\sqrt{T_0/6}$.  We apply Corollary~\ref{corrlqlp} to the rescaled solutions $a^{\eps}$.

Each component  $a^{\eps}_i\in L^q(-3,0;L^p (B(0,3)))$ for  $p\in(1,\frac{N}{N-2})$ and $q$ such that $\frac{1}{p}=1-\frac{2}{qN}$; Moreover, for each $1\leq i\leq P$
$$
\|a^{\eps}_i\|_{L^q(-3, 0;L^p(B(0,3)))}\leq C \|\nabla_x\sqrt { a^{\eps}_i}\|_{ L^2((-3, 0)\times R^N)}^{\frac{2(p-1)}{p}}\|\rho^{\eps}\|_{L^\infty(-3,0;\mathcal{L}_w(\R^N))}^{\frac{1}{p}},
$$

\begin{equation}\label{normeps2}
\|\nabla_x\sqrt{a_i^{\eps}}\|^{2}_{L^2((-3,0)\times\R^N)}\leq\eps^{\frac{2}{\nu-1}-N}\|\nabla_x\sqrt{a_i}\|^{2}_{L^2((0,T_0)\times\R^N)}.
\end{equation}
Hence, together with (\ref{final control weak norm}), we obtain:

$$\sum_{i=1}^{P} \|a^{\eps}_i\|_{L^q(-3,0;L^p (B(0,3)))}\leq C(M_0,T_0) \eps^{\alpha(p)} $$
\noindent where $\alpha(p):=\frac{p-1}{p}(\frac{2}{\nu-1}-2)+\frac{1}{p}(\frac{2}{\nu-1}-N)$ and $ C(M_0,T_0)$ is a constant depending on the quantity $M_0$ defined by~(\ref{M_0}) and $T_0$.\\

\noindent Since $\nu<2$, the factor~($\frac{2}{\nu-1}-2$) is positive. Thus, there exists $\bar p\in(1,\frac{N}{N-2})$ such that $\bar \alpha :=\alpha(\bar p)>0$ and
$$\sum_{i=1}^{P} \|a^{\eps}_i\|_{L^{\bar q}((-3,0)\times B(0,3))}\leq C(M_0,T_0) \eps^{\bar\alpha} $$
 with $\bar q=\frac{2}{N}\big(\frac{\bar p}{\bar p-1}\big)$. This last inequality implies that it is enough to choose $\eps_0\leq \inf\left\{\sqrt{\frac{T_0}{6}},\left(\frac{\delta^{\star}}{C(M_0,T_0)}\right)^{\frac{1}{\bar \alpha}}\right\}$ to conclude the statement for all $0<\eps\leq\eps_0$.

\endproof
\vskip1cm
Next, we apply Proposition~\ref{lemm_de giorgi} to the rescaled functions $a^{\eps}$.
\vskip1cm
\begin{corr}\label{unibound} Let $\delta^{\star}=\delta^{\star}(\bar q)$, $\bar\alpha$ and $C(M_0,T_0)$ as in Lemma~\ref{lemm_de giorgi2}. Then, for all $\frac{T_0}{2}<T<T_0$, we have
$$\sum_{i=1}^{P} ||a_i(T,\cdot)||_{L^{\infty}(\R^N)}\leq \left(\inf\left\{\sqrt{\frac{T_0}{6}},\left(\frac{\delta^{\star}}{C(M_0,T_0)}\right)^{\frac{1}{\bar \alpha}}\right\}\right)^{-\frac{2}{\nu-1}}.$$
\end{corr}
\begproof
Let $x_0$ be any point in $\R^N$,  $\frac{T_0}{2}<T<T_0$. Let  $a^{\eps_0}$ be the rescaled function defined in~(\ref{eq_a epsilon}) with $\eps_0=\inf\left\{\sqrt{\frac{T_0}{6}},\left(\frac{\delta^{\star}}{C(M_0,T_0)}\right)^{\frac{1}{\bar\alpha}}\right\}$. Then, by Lemma~\ref{lemm_de giorgi2}, we can apply Proposition~\ref{lemm_de giorgi} to $a^{\eps_0}$ to obtain:

$$0\leq \sum_{i=1}^{P}  a_i(x_0,T)=\eps_0^{\frac{-2}{\nu-1}}\sum_{i=1}^{P} a_i^{\eps_0}(0,0)\leq \eps_0^{\frac{-2}{\nu-1}}P.$$

\endproof

{ \it Finally, we can now prove Theorem~1.}\\

Assume by contradiction that the maximal time of existence of a solution $a=(a_1,\cdot\cdot\cdot, a_P)$, $T_0$ is finite. Then, by Corollary~\ref{unibound}, each $a_i(t,\cdot)$, for $1\leq i\leq P$ would be uniformly bounded for all $\frac{T_0}{2}t<T_0$. In particular, it would be bounded at $t=T_0$. Hence, by standard arguments this would imply that $a\in C^{\infty}$ at $t=T_0$. This concludes the proof because it would negate the fact that $T_0$ is the maximal time of existence of a smooth solution $a=(a_1,\cdot\cdot\cdot, a_P)$.

\appendix \label{appendix}
\section{Proof of Theorem \ref{theo_max_principle}}
This appendix is dedicated to the proof of the maximum principle in the case $P=2$. This is a very standard proof. We include it here since we did not find the result in the literature.
The system is equivalent to the pair of equations
\begin{equation}\label{eq_system1}
\left\{\begin{array}{l}
\displaystyle{\dt a_1-d_1\Delta a_1=Q(a), \qquad t>0, \ x\in
\R^N,}\\[0.3cm]
\displaystyle{\dt a_2-d_2\Delta a_2=-Q(a), \qquad t>0, \ x\in
\R^N,}\\[0.3cm]
\displaystyle{a(0,x)=a^0(x), \qquad x\in \R^N.}
\end{array}
\right.
\end{equation}
The main remark is that, under Hypothesis (\ref{Hypothesis_entropy}), we have
\begin{equation}\label{signe de Q}
Q(a_1,a_2)(a_1-a_2)\leq 0, \qquad \mathrm{for\ any\ } a_1,a_2\geq0.
\end{equation}
Let $(a_1,a_2)$ be a smooth solution ($C^2$) to (\ref{eq_system1}) on $([0,T]\times\R^N)$, decaying to $0$ when $|x|\to\infty$.
For any $\eps>0$, we define
\begin{eqnarray*}
&&a_1^\eps(t,x)=a_1(t,x)-\eps t,\\
&&a_2^\eps(t,x)=a_2(t,x)-\eps t.
\end{eqnarray*}
Assume that $\sup(a_1^\eps,a_2^\eps)$ attains a local maximum at a point $(t_\eps,x_\eps)\in (0,T]\times\R^N$. Let say that its value is $a_1^\eps(t_\eps,x_\eps)$. (The proof is similar for the other case.) Especially, this is a local maximum for $a_1^\eps$. Note that
$$
a_1^\eps-a_2^\eps=a_1-a_2.
$$
Therefore, from (\ref{signe de Q}), we have
$$
Q(a)(t_\eps,x_\eps)\leq0.
$$
From the first equation of (\ref{eq_system1}) we have
$$
(\dt a^\eps_1-d_1\Delta a_1^\eps)(t_\eps,x_\eps)=Q(a)(t_\eps,x_\eps)-\eps\leq-\eps<0.
$$
This contradicts the fact that $a_1^\eps$ attains a local maximum at $(t_\eps,x_\eps)$.
Hence
$$
\sup_{(t,x)\in(0,T]\times\R^N}(a_1^\eps(t,x),a_2^\eps(t,x))\leq \sup_{x\in\R^N}(a_1^\eps(t=0,x),a_2^\eps(t=0,x))=\sup_{x\in\R^N}(a_1^0(x),a_2^0(x)).
$$
Passing to the limit as $\eps\to0$ gives the result.

\end{document}